\documentclass{article}
\usepackage{amssymb,amsmath}
\usepackage{graphicx}

\def\qed{\hfill {\hbox{${\vcenter{\vbox{               
   \hrule height 0.4pt\hbox{\vrule width 0.4pt height 6pt
   \kern5pt\vrule width 0.4pt}\hrule height 0.4pt}}}$}}}

\def\tr{\triangleright}
\def\bar{\overline}

\newtheorem{theorem}{Theorem}
\newtheorem{definition}{Definition}
\newtheorem{lemma}[theorem]{Lemma}
\newtheorem{proposition}[theorem]{Proposition}
\newtheorem{corollary}[theorem]{Corollary}
\newtheorem{example}{Example}

\newtheorem{remark}{Remark}

\newenvironment{proof}[1][Proof]{\smallskip\noindent{\bf #1.}\quad}%
{\qed\par\medskip}

\date{}

\title{\Large \textbf{Virtual Yang-Baxter cocycle invariants}}

\author{
\begin{tabular}{c} Jose Ceniceros  \\ 
\small Whittier College \\
\small 13406 Philadelphia \\
\small P.O. Box 634 \\
\small Whittier, CA 90608-0634 \\
\texttt{jcenicer@poets.whittier.edu}\end{tabular}
\and
\begin{tabular}{c} Sam Nelson \\ 
\small Department of Mathematics \\
\small Pomona College \\
\small 610 North College Avenue \\ 
\small Claremont, CA 91711 \\
\texttt{knots@esotericka.org}\end{tabular}
}

\begin{document}
\maketitle

\begin{abstract}
We extend the Yang-Baxter cocycle invariants for virtual knots by augmenting
Yang-Baxter 2-cocycles with cocycles from a cohomology theory associated
to a virtual biquandle structure. These invariants coincide with the classical
Yang-Baxter cocycle invariants for classical knots but provide extra 
information about virtual knots and links. In particular, they provide a 
method for detecting non-classicality of virtual knots and links. 

\end{abstract}

\textsc{Keywords:} Finite biquandles, virtual biquandles, virtual link 
invariants, Yang-Baxter cohomology

\textsc{2000 MSC:} 57M27, 18G60

\section{\large \textbf{Introduction}}

In \cite{CYB}, quandle 2-cocycle invariants are generalized to the 
biquandle case. These are link invariants associated to finite biquandles
and 2-cocycles in their Yang-Baxter cohomology. In \cite{KM}, biquandles 
are generalized to include an operation at virtual crossings, a structure 
called a \textit{virtual biquandle}.   

A \textit{virtual link} is an equivalence class of link diagrams with 
classical crossings and virtual crossings under virtual Reidemeister moves.
Virtual crossings are crossings which arise from genus in the surface on 
which the knot diagram is drawn, as distinguished from ordinary 
\textit{classical} crossings.

Biquandles are algebras with axioms derived from the Reidemeister moves.
Associated to a finite biquandle $T$ is the \textit{counting invariant} 
$|\mathrm{Hom}(B(L),T)|$, the number of colorings of a link diagram by $T$. 
The Yang-Baxter 2-cocycle invariants are jazzed-up versions of this counting 
invariant where we add up the values of a function $\phi$ at every crossing 
in each biquandle-colored diagram to obtain a \textit{Boltzmann weight} for 
that coloring. Such a function $\phi$ defines an invariant provided 
$\phi$ is a Yang-Baxter 2-cocycle.

Including a virtual operation in a finite biquandle structure strengthens 
the counting invariants for virtual links. In this paper we describe a
new method for strengthening and extracting more information from the
virtual biquandle counting invariants using Yang-Baxter cocycles together 
with cocycles from a cohomology theory associated to the virtual operation. 
In particular, the new invariants can detect non-classicality in some 
virtual links.

The paper is organized as follows. In section 2 we recall the basics of 
virtual knot theory; in section 3 we recall virtual quandles and virtual 
biquandles. In section 4 we recall Yang-Baxter cohomology. In section
5 we introduce $S$-homology and our virtual Yang-Baxter cocycle invariants.
In section 6 we see some examples of the new invariants, and in section
7 we collect some questions for future research.

\section{\large \textbf{Classical and virtual knots and links}}

A \textit{classical oriented link diagram} is a planar 4-valent directed
graph with two edges directed in and two directed out at every vertex. The 
vertices are regarded as crossings, so one inbound/outbound pair forms an 
over-crossing strand and the other pair forms an under-crossing strand; 
we indicate which is which by drawing the under-crossing strand broken. 
It is well-known that 
two classical oriented link diagrams represent ambient isotopic links iff 
they are related by a finite sequence of oriented Reidemeister moves. Here 
we depict the unoriented moves:

\[\includegraphics{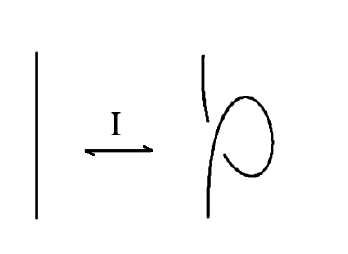}\quad
\includegraphics{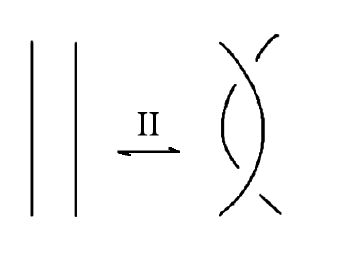}\quad
\includegraphics{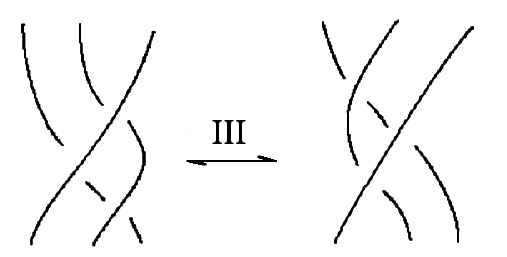}\quad
\]

A \textit{virtual oriented link diagram} is an oriented link diagram in which 
we drop the planarity requirement. We may draw such a diagram on an 
orientable surface $\Sigma$ with genus $g\ge 0$, and the represented 
\textit{virtual link} is really an ordinary link whose ambient space is 
$\Sigma\times [0,1]$. Classical links form a subset of virtual links. 
To draw an oriented virtual link on planar paper, we must distinguish 
crossings arising from genus in $\Sigma$ from ordinary \textit{classical} 
crossings (decorated vertices in the non-planar graph); the former are called 
\textit{virtual crossings}, represented as circled self-intersections with 
no over/under sense. 
Virtual links were introduced in \cite{K} and have been studied in many 
recent works.

\[\includegraphics{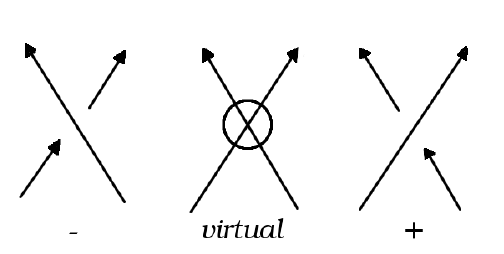}\quad\quad \quad
\includegraphics{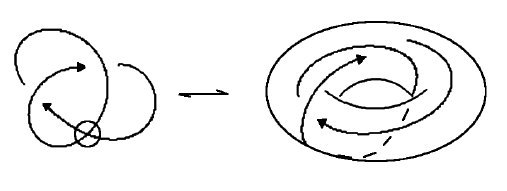}\]

We may also represent an oriented link diagram using a \textit{signed Gauss 
code} in which every classical crossing is assigned a label and a plus or 
minus sign according to its local writhe number or \textit{crossing sign}. 
A signed Gauss code is then the result of traveling around 
each component of the oriented link diagram in the direction indicated by the 
orientation and noting the order in which each overcrossing and undercrossing 
is encountered. 

\begin{example}
The figure eight knot and its signed Gauss code.
\[\raisebox{-0.5in}{\includegraphics{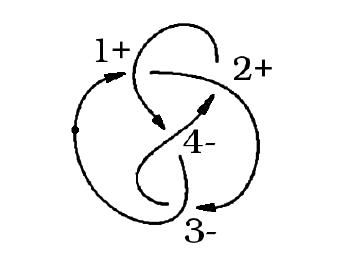}} 
\quad U1^+O2^+U3^-O4^-U2^+O1^+U4^-O3^-\]
\end{example}

Indeed, signed Gauss codes provide one motivation for studying virtual
links -- while not every Gauss code corresponds to a planar link diagram,
every signed Gauss code can be interpreted as a virtual link diagram by 
simply inserting virtual crossings as needed while drawing the link. A 
\textit{virtual link} is then an equivalence class of oriented virtual link
diagrams or, equivalently, signed Gauss codes, under the equivalence relation 
generated by the three Reidemeister moves.

Since virtual crossings are not included in ordinary signed Gauss codes,
any two oriented virtual link diagrams with the same signed Gauss code
are equivalent. Such diagrams are related by the \textit{detour move}, which
can be broken down into four \textit{virtual Reidemeister moves}.

\[\includegraphics{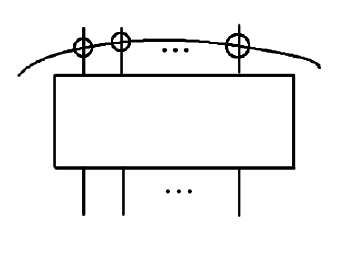} \quad
\raisebox{0.5in}{$\stackrel{\mathrm{Detour}}{\longleftrightarrow}$}
\quad \includegraphics{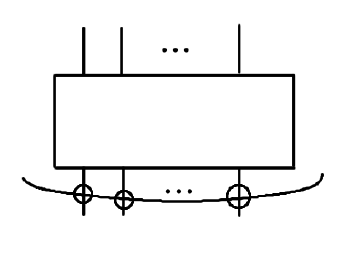}
\quad\quad\quad \includegraphics{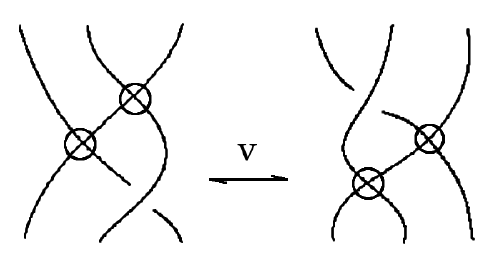}
\]

\[
\scalebox{0.95}{\includegraphics{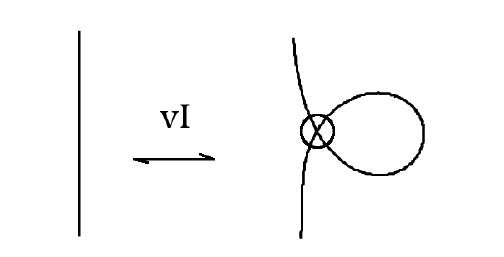}} \quad
\scalebox{0.95}{\includegraphics{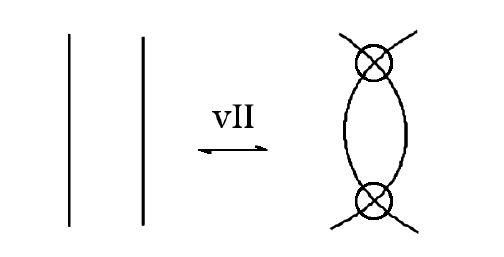}} \quad
\scalebox{0.95}{\includegraphics{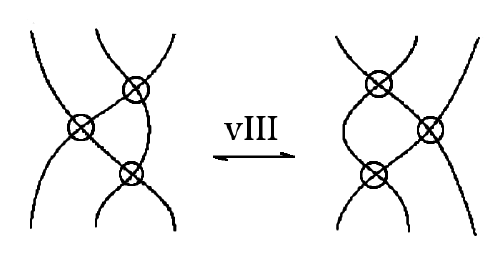}} \quad
\]

We may then consider virtual knots to be equivalence classes of planar 
virtual knot diagrams under the equivalence relation generated by the three
Reidemeister moves and the detour move, or equivalently the three
Reidemeister moves and the four virtual moves $v, vI, vII$ and $vIII$. Note
that while move $v$ permits a strand with only virtual crossings to be moved
past a classical crossing, we cannot move a strand with classical crossings
over or under a virtual crossing, since this changes the Gauss code and,
generally, the virtual knot type.

Virtual crossings are not usually included in a signed Gauss code, since
in one sense the whole point of virtual crossings is to avoid testing
codes for planarity. However, to compute the new invariants we will define 
in section \ref{inv}, we will need to specify the location of our virtual 
crossings. Thus, we have

\begin{definition}
\textup{A \textit{virtual signed Gauss code} is a signed Gauss code in 
which we specify the virtual crossings in addition to the classical crossings,
with ``R'' and ``L'' labels indicating whether we are entering the
virtual crossing from the right or left input.}
\end{definition}

\begin{example}
\textup{A virtual knot diagram and its virtual signed Gauss code.}
\[\includegraphics{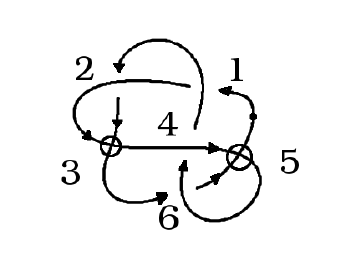}\quad 
\raisebox{0.5in}{$U1^+O2^+R3\ O4^+L5\ O6^-U4^+O1^+U2^+L3\ U6^-R5$}
\]
\end{example}

The following observation will be useful.

\begin{lemma}\label{lem:v}
The classical and virtual II moves, together with one oriented v
or vIII move, imply the other oriented v and vIII moves. That is, 
we can reverse the direction of any strand in a type v or vIII move
using a sequence of II and vII moves. 
\end{lemma}

\begin{proof}
We illustrate the case of reversing the virtual strand in a v move using 
virtual II moves. The other cases are similar.
\[\includegraphics{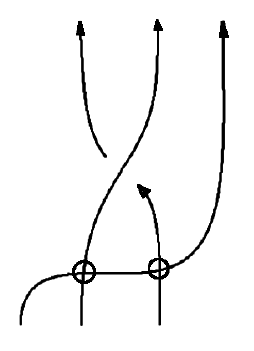} \raisebox{0.75in}{$\leftrightarrow$} 
\includegraphics{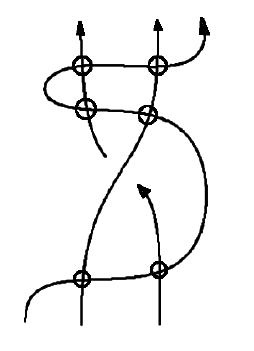} \raisebox{0.75in}{$\leftrightarrow$} 
\includegraphics{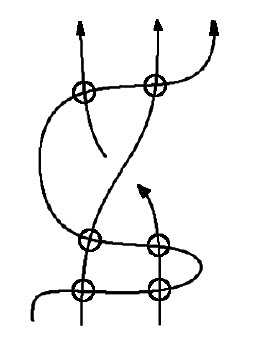} \raisebox{0.75in}{$\leftrightarrow$} 
\includegraphics{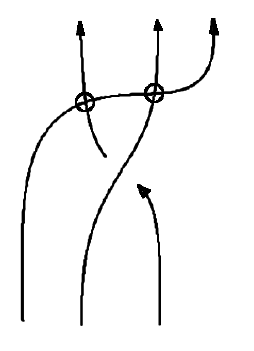}
\]
\end{proof}

\section{\large \textbf{Virtual quandles and biquandles}}

The combinatorial approach to knot theory described in the previous 
section has a corresponding algebraic approach, which we now describe.

\begin{definition}
\textup{
A \textit{biquandle} is a set $B$ with four binary operations denoted
$(a,b)\mapsto a^b, a^{\bar b}, a_b, a_{\bar b}$ satisfying the following 
axioms:
\begin{list}{}{}
\item[1.]{For every $a,b\in B$ we have
\[a=(a^b)^{(\bar{b_a})}, \quad b=(b_a)_{(\bar{a^b})}, 
\quad a=(a^{\bar{b}})^{(b_{\bar{a}})},
\quad \mathrm{and}
\quad b=(b_{\bar{a}})_{(a^{\bar{b}})},\]}
\item[2.]{for every $a,b\in B$ there exist unique $x,y\in B$ such that
\[ x=a^{(b_{\bar{x}})},\quad  a=x^{\bar{b}},\quad  b=(b_{\bar{x}})_a,\quad
y=a^{(\bar{b_y})},\quad a=y^b,\quad  \mathrm{and} \quad b=(b_y)_{\bar{a}},\]}
\item[3]{for every $a,b,c\in B$ we have
\[(a^b)^c=(a^{(c_b)})^{(b^c)}, \quad (c_b)_a=(c_{(a_b)})_{(b_a)}, \quad 
(b_a)^{(c_{(a^b)})}=(b^c)_{(a^{(c_b)})},\]
\[(a^{\bar{b}})^{\bar{c}}=(a^{\bar{(c_{\bar{b}})}})^{\bar{(b^{\bar{c}})}}, 
\quad 
(c_{\bar{b}})_{\bar{a}}=(c_{\bar{(a_{\bar{b}})}})_{\bar{(b_{\bar{a}})}}, 
\quad 
\mathrm{and}\quad 
(b_{\bar{a}})^{\bar{(c_{\bar{(a^{\bar{b}})}})}}
=(b^{\bar{c}})_{\bar{(a^{\bar{(c_{\bar{b}})}})}},\]}
\item[4.]{for every $a\in B$ there exist unique $x,y\in B$ such that
\[x=a_x, \quad a=x^a, \quad y=a^{\bar{y}}\quad \mathrm{and} \quad 
a=y_{\bar{a}}.\]}
\end{list}
A biquandle in which $a_b=a_{\bar{b}}=a$ for all $a,b\in B$ is a 
\textit{quandle}; if $B$ is a quandle we may denote $a^b=a\tr b$ and
$a^{\bar{b}}=a\tr^{-1} b$.}
\end{definition}

The biquandle axioms are obtained by associating elements of $B$ to
semiarcs (i.e., edges in the link diagram considered as a 4-valent graph) 
in an oriented link diagram and letting these elements act on each other 
at crossings as pictured. The notation is chosen so that $a^b$ suggests
``$a$ under $b$'' and $b_a$ suggests ``$b$ over $a$.''

\[\includegraphics{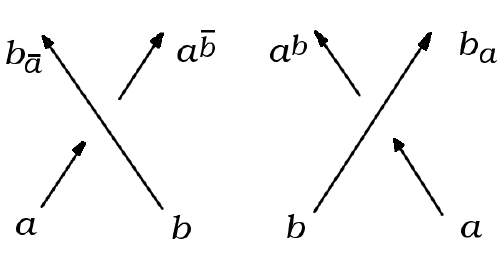}\]

The set of biquandle axioms is then the result of checking the conditions 
imposed by a minimal set of oriented Reidemeister moves. For example, the 
type II Reidemeister move with both strands oriented in the same direction 
(the \textit{direct} II-move) implies axiom (1) above:
\[\includegraphics{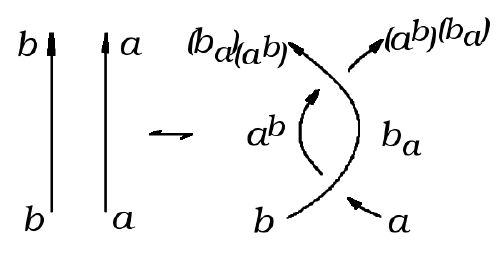}\quad \quad 
\includegraphics{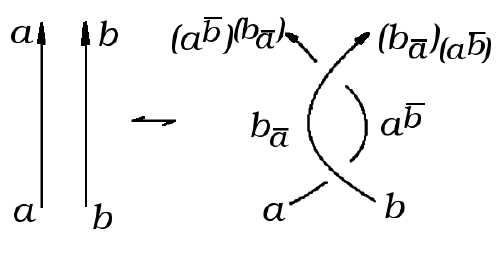}\]
See \cite{KR} for more.

Biquandles are useful for defining invariants of knots and links. 
Specifically, for every knot or link diagram there is an associated
\textit{knot biquandle} defined as the quotient of the free biquandle
on the set of semiarcs in the diagram by the equivalence relation 
generated by the relations determined by the crossings in the diagram.

Given a finite biquandle $T$, we can obtain an easily computable link
invariant by computing the cardinality $|\mathrm{Hom}(B(L),T)|$ of the 
set of biquandle homomorphisms from the knot biquandle $B(L)$ of the
link $L$ into $T$. Each such homomorphism can be pictured as a ``coloring''
of $L$ which associates an element of $T$ to each semiarc in $L$ such that
the biquandle operations in $T$ are compatible with the relations determined
by the crossings in $L$.

Thus, for every finite biquandle $T$ there is a counting invariant 
$|\mathrm{Hom}(B(L),T)|$. Several recent papers describe examples of finite
biquandles (see \cite{CS}, \cite{FJK} and \cite{NR} for example). One nice
example is the family of \textit{Alexander biquandles}, where $T$ is a 
module over the ring $\mathbb{Z}[s^{\pm 1},t^{\pm 1}]$ of two-variable
Laurent polynomials with integer coefficients. Such a $T$ is a biquandle with 
operations
\[ a^b=ta+(1-st)b, \ a^{\bar{b}}=t^{-1} a + (1-s^{-1}t^{-1}) b, \
a_b=sa, \ a_{\bar{b}}=s^{-1}a.\]
For example, for any $n\ge 2$ we can get an Alexander biquandle of
cardinality $n$ by taking $T=\mathbb{Z}_n$ and choosing $s,t$ invertible 
elements of $\mathbb{Z}_n$.

\begin{example}\label{tref}
\textup{The trefoil knot $3_1$ has three colorings by the Alexander 
biquandle $\mathbb{Z}_3$ with $s=t=2$:}
\[\includegraphics{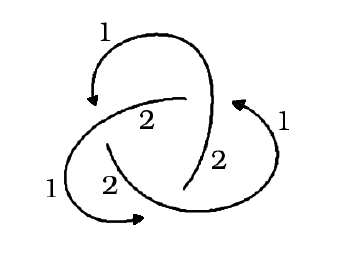}\quad 
\includegraphics{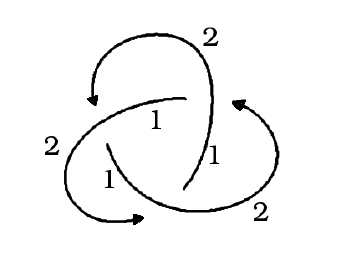}\quad 
\includegraphics{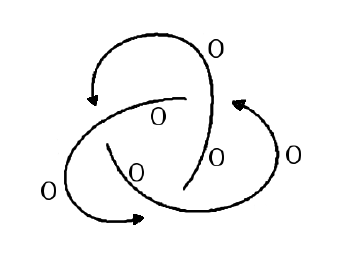}
\]
\textup{Hence, we have $|\mathrm{Hom}(B(3_1),T)|=3.$}
\end{example}

We can specify a finite biquandle $T=\{x_1,\dots, x_n\}$ by listing its 
operation tables in the form of a block matrix 
$M_T=\left[\begin{array}{c|c}
i^{\bar{j}} & i^j \\ \hline
i_{\bar{j}} & i_j 
\end{array}\right]$; for example the biquandle $T=\mathbb{Z}_3$ with 
$s=t=2$ from example \ref{tref} has biquandle matrix
\[M_T=\left[\begin{array}{ccc|ccc}
2 & 2 & 2 & 2 & 2 & 2 \\
1 & 1 & 1 & 1 & 1 & 1 \\
3 & 3 & 3 & 3 & 3 & 3 \\ \hline
2 & 2 & 2 & 2 & 2 & 2 \\
1 & 1 & 1 & 1 & 1 & 1 \\
3 & 3 & 3 & 3 & 3 & 3 \\
\end{array}\right]\]
where we use $x_1=[1], x_2=[2]$ and $x_3=[0]$. This notation has the 
advantage of letting us compute the counting invariant even in the
absence of a convenient algebraic description of the structure of $T$.
See \cite{NV} and \cite{CN} for more.

For virtual knots, we can simply ignore the presence of virtual crossings and
compute the counting invariant with respect to a finite biquandle $T$, since
the Reidemeister moves are local and hence agnostic of the genus of the 
surface on which $L$ is drawn.

An alternative idea is to consider virtual crossings in our algebraic
structure, dividing our virtual link diagrams at both classical and virtual
crossings and introducing an additional operation at the virtual crossings.
This approach was taken in \cite{KM}; the resulting algebraic structure is 
called a \textit{virtual biquandle}.

\begin{definition}\textup{
Let $B$ be a biquandle and let $S:B\to B$ be an automorphism of $B$,
that is, $S$ is a bijection such that for all $a,b\in B$}
\[S(a^b)=S(a)^{S(b)}, \ S(a_b)=S(a)_{S(b)}, 
\ S(a^{\bar{b}})=S(a)^{\bar{S(b)}}, \ \mathrm{and} \ 
\ S(a_{\bar{b}})=S(a)_{\bar{S(b)}}.
\]
\textup{Then the pair $(B,S)$ is a \textit{virtual biquandle}. 
If $a_b=a_{\bar{b}}=a$ for all $a\in B$, then $B$ is a quandle and $(B,S)$
is a \textit{virtual quandle}.}
\end{definition}

At a virtual crossing, the semiarc entering from the right picks up an
$S$ while the semiarc entering from the left picks up an $S^{-1}$:
\[\includegraphics{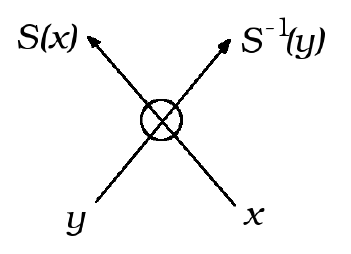}\]

The pure virtual moves $vI$ and $vII$ require that $S$ is bijective, 
move $vIII$ imposes no conditions, and the oriented $v$ moves
require that $S$ is an automorphism of $B$; see \cite{KM} (or, 
check it as an exercise!).

We may regard ordinary biquandles as the special case $S=Id_B$, so the 
theory of biquandles is a subset of the theory of virtual biquandles.
Given a virtual link diagram $L$, one defines the \textit{knot virtual 
biquandle} $VB(L)$ of $L$ as the quotient of the free virtual biquandle 
on the set of semiarcs in the virtual link diagram $L$ by the equivalence
relation generated by the relations at the crossings in $L$. As before, 
we can regard a homomorphism $f:VB(L)\to (T,S)$ from the knot virtual 
biquandle of $L$ to a finite virtual biquandle $(T,S)$ as a coloring of $L$,
i.e., an assignment of elements of $T$ to semiarcs in $L$ (divided at 
both classical and virtual crossings) such that the virtual biquandle
operations are satisfied at every crossing.


\begin{proposition}
Let $B$ be a biquandle and $S,S'\in \mathrm{Aut}(B)$. Then the two virtual 
biquandle structures $(B,S)$ and $(B,S')$ are isomorphic iff $S$ is 
conjugate to $S'$. 
\end{proposition}

\begin{proof}
An isomorphism of virtual biquandles $\phi:(B,S)\to (B,S')$
is an automorphism $\phi:B\to B$ such that for all $x\in B$ we have 
$\phi(S(x))=S'(\phi(x))$, that is, $\phi S= S' \phi$
or $S=\phi^{-1} S' \phi.$ Conversely, if $S=\phi^{-1}S'\phi$ in
$\mathrm{Aut}(B)$, then $\phi(S(x))=S'(\phi(x))$ for all $x\in B$ and
$\phi$ is an isomorphism from $(B,S)$ to $(B,S').$
\end{proof}

\begin{corollary}
Given a biquandle $B$, we may list all virtual biquandle structures on $B$
by listing one representative from each conjugacy class of $\mathrm{Aut}(B).$
\end{corollary}

\begin{example}
\textup{The Alexander biquandle $T=\mathbb{Z}_3$ with $s=1$, $t=2$ has
biquandle matrix}
\[M_T=\left[\begin{array}{ccc|ccc}
1 & 3 & 2 & 1 & 3 & 2 \\
3 & 2 & 1 & 3 & 2 & 1 \\
2 & 1 & 3 & 2 & 1 & 3 \\ \hline
1 & 1 & 1 & 1 & 1 & 1 \\
2 & 2 & 2 & 2 & 2 & 2 \\
3 & 3 & 3 & 3 & 3 & 3 
\end{array}\right]\]
\textup{In fact, this is an Alexander quandle, $\mathbb{Z}[t^{\pm 1}]/(t-2)$.
It has automorphism group $S_3$, with conjugacy classes represented by
$(1)$, $(23)$ and $(123)$. Hence, up to isomorphism there are three virtual
quandle structures on $T$.}
\end{example}

If a virtual diagram $L$ represents a classical knot or link, then the 
counting invariants with respect to any two virtual biquandle structures 
on a given biquandle must be equal since both can be obtained from a diagram
equivalent to $L$ with no virtual crossings. Thus we have:

\begin{proposition}
If $|\mathrm{Hom}(VB(L),(T,S))|\ne |\mathrm{Hom}(VB(L),(T,S'))|$ for any
two virtual biquandle structures $(T,S)$ and $(T,S')$ defined on the same 
biquandle $T$, then $L$ is non-classical.
\end{proposition}

\begin{example}
\textup{The virtual quandle structures on the Alexander quandle 
$T=\mathbb{Z}[t^{\pm 1}]/(t-2)$ show that the virtual trefoil is 
non-classical, since $(T,(1))$ yields a counting invariant value of 3
while $(T,(123))$ has a counting invariant value of 0.}
\[\includegraphics{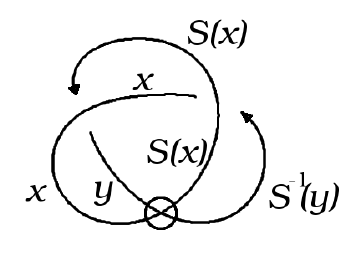} \quad\raisebox{0.5in}{$
\begin{array}{rcl}
2S^{-1}(y)+2S(x) & = & x \\
2S(x)+2x & = & y 
\end{array}$}\]
\textup{Here we have $x^y=x^{\bar{y}}=2x+2y \ (\mathrm{mod} \ 3)$ and
$x_y=x_{\bar{y}}=x$. Thus, colorings are solutions in $\mathbb{Z}_3$ 
to the above system of equations. 
If $S(x)=x$, the system reduces to $\begin{array}{rcl}
y & = & x \\
x & = & y
\end{array}$
and there is one solution for every element of $T$. On the other
hand, the permutation $S=(123)$ translates as $S(x)=x+1$, and our system
becomes 
\[\begin{array}{rcl}
2y-2+2x+2 & = & x \\
2x+2+2x & = & y
\end{array}\quad \mathrm{ which \ reduces \ to} \quad
\begin{array}{rcl}
y & = & x \\
x+2 & = & y
\end{array}\]
which is inconsistent.}
\end{example}

The counting invariants associated to virtual biquandles provide a 
sharpening and strengthening of the biquandle counting invariants
for virtual knots and links, while they coincide with the usual biquandle
counting invariants for classical knots and links. In the next section,
we will see a method for further improving and strengthening these
invariants.

\section{\large \textbf{Yang-Baxter Cohomology}}

The counting invariant $|\mathrm{Hom}(B(L),T)|$ associated to a finite
biquandle $T$ is a computable and useful invariant of knots and links. 
However, a set is more than a mere cardinality, and it is natural to seek
ways of extracting more information from the set $\mathrm{Hom}(B(L),T)$.

One very successful idea, with origins in statistical mechanics, is to 
associate to every crossing a function $\phi$ with values in some abelian 
group; then for a given element $f\in \mathrm{Hom}(B(L),T)$, we can sum 
our values of $\phi$ over the crossings in the diagram colored by $f$ to 
obtain a \textit{Boltzmann weight} for $f$. The set-with-multiplicities 
of all such weights for all $f\in \mathrm{Hom}(B(L),T)$ is then an 
invariant of knots and links provided the function $\phi$ is chosen such 
that the Boltzmann weight for each colored diagram is unchanged by
Reidemeister moves.

The condition on $\phi$ required to make this happen turns out to be
expressible in terms of a cohomology theory associated to the coloring
quandle or biquandle $T$. Many variations of this idea exist, including
the quandle homology described in \cite{CJKLS}, twisted
quandle homology described in \cite{TWIST}, quandle homology with 
coefficients in a quandle module (\cite{AG}, \cite{J}), the biquandle version
known as Yang-Baxter homology (\cite{CYB}), and a unifying approach to
the homology of quandles and Lie algebras (\cite{CC}).

Since our focus is on virtual biquandles, we will use the Yang-Baxter 
homology approach from \cite{CYB}. Moreover, we will choose our coefficient 
ring to be $\mathbb{Q}$ for simplicity of computation. Then if $B$ is a 
biquandle, a function $\phi:B\times B\to \mathbb{Q}$ is a 
\textit{Yang-Baxter 2-cocycle} if $\phi$ satisfies
\[\phi(a,b)+\phi(a^b,c)+\phi(b_a, c_{a^b}) =
\phi(b,c)+\phi(a,c_b)+\phi(a^{c_b},b^c) \]
for all $a,b,c\in B$. If $\phi(x,a)=\phi(a,y)=0$ for all $a\in B$ where
$x,y$ are the elements of $B$ from biquandle axiom (4), then $\phi$ is
a \textit{reduced Yang-Baxter 2-cocycle}. We will refer to such
cocycles as \textit{classical cocycles}.

As implied by the term ``2-cocycle'', these functions $\phi$ are elements 
of the second cohomology in a cohomology theory, described in \cite{CYB}.
The chain groups are generated by ordered $n$-tuples of biquandle elements,
and the boundary map is obtained by interpreting an $n$-tuple as colors
on a preferred path in the $n$-cube from $(0,0,\dots, 0)$ to 
$(1,1,\dots, 1)$. Such a coloring then extends in a unique way to the whole 
$n$-cube graph; one then takes the usual boundary of the $n$-cube and 
interprets the faces as the $(n-1)$-tuples corresponding to their 
preferred paths. It is perhaps more enlightening, however, to consider
the Reidemeister III move:

\[\includegraphics{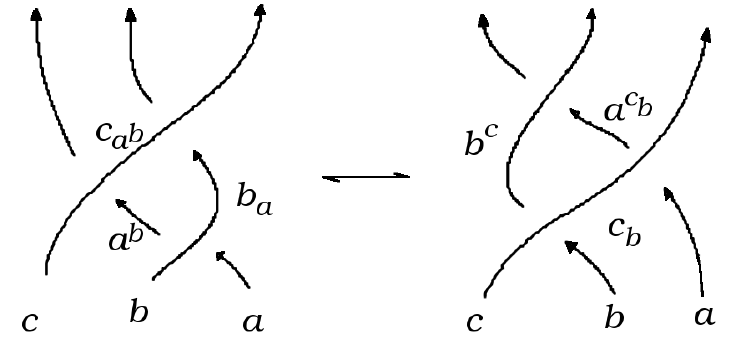} \quad \quad\quad 
\raisebox{0.15in}{ \includegraphics{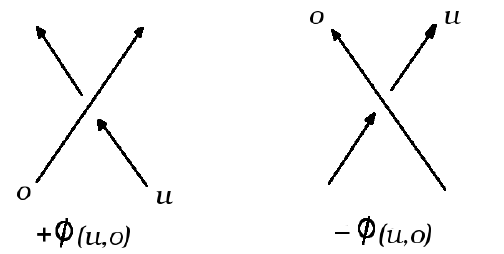}}\]

Here the rule is that every positive crossing contributes $\phi(u,o)$
where $u$ is the color on the undercrossing inbound semiarc and $o$
is the color on the overcrossing inbound semiarc, while negative 
crossings contribute $-\phi(u,o)$ where $u$ is the undercrossing outbound
semiarc, and $o$ is the overcrossing outbound semiarc.

This rule guarantees that the Boltzmann weights are unchanged by 
Reidemeister II moves; the cocycle condition takes care of the type III
moves and the reduced condition takes care of the type I moves. 
The classical Yang-Baxter 2-cocycle invariant is then
\[\Phi_{YB}(L) = \sum_{f\in \mathrm{Hom}(B(L),T)} t^{BW(f)}\]
where $BW(f)=\sum_{\mathrm{crossings}} \pm \phi(u,o)$ is the Boltzmann 
weight of the homomorphism $f$.

Finally, it is observed in \cite{CYB} that if 
$f=\delta^1_{YB}(g)=g(\partial^{YB}_2)$ is a Yang-Baxter coboundary
then at every crossing $f$ effectively counts $g(ui)+g(oi)-g(uo)-g(oo)$ 
where $ui, \ oi, \ uo$ and $oo$ are the inbound undercrossing, inbound 
overcrossing, outbound undercrossing, and outbound overcrossing semiarcs 
respectively. Then since every semiarc appears in the sum once as an inbound 
and once as an outbound semiarc, the total Boltzmann weight of a coboundary 
is zero and cohomologous cocycles define the same invariant. 

\section{\large \textbf{$S$-homology}}

For virtual biquandles $(B,S)$, we now define a second cohomology theory which
we will use to augment the classical theory for defining invariants of 
virtual knots and links. First, we need a rule for contributions to the
Boltzmann weight from each virtual crossing in a virtual biquandle colored 
diagram. Let $v:B\times B\to \mathbb{Q}$; then at each virtual crossing, 
we will count $v(x,y)-v(S^{-1}(y),S(x))$ where $x$ is the biquandle color 
on the right-hand input semiarc, $y$ is the color on the left-hand input 
semiarc, $S^{-1}(y)$ is the output right-hand semiarc and $S(x)$ is the 
color on left-hand output semiarc.

\[\includegraphics{jc-sn-24.pdf} \quad \quad 
\raisebox{0.5in}{$v(x,y)-v(S^{-1}(y),S(x))$}\]

\begin{proposition}
The contribution to the Boltzmann weight of a colored virtual link diagram
with the above rule is unchanged by vI and vII moves.
\end{proposition}

\begin{proof}
Move $vI$:
\[\includegraphics{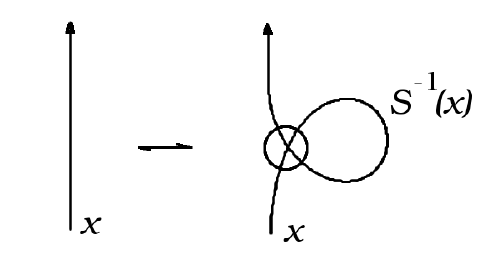}\hskip 1in \includegraphics{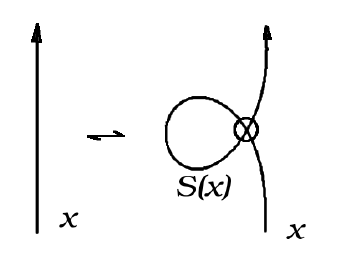}\]
In the first move, the contribution to the Boltzmann weight on the left is 0,
and on the right, the contribution is $v(S^{-1}(x),x)-v(S^{-1}(x),x)=0.$ 
In the second move, on the left again we have a contribution of 0 and on 
the right we have $v(x,S(x))-v(x,S(x))=0$.

Move $vII$: There are two cases, one with the two strands oriented in the 
same direction and one with the two strand oriented in opposite directions.
\[\includegraphics{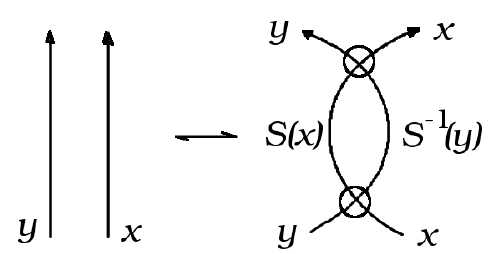} \hskip 1in  \includegraphics{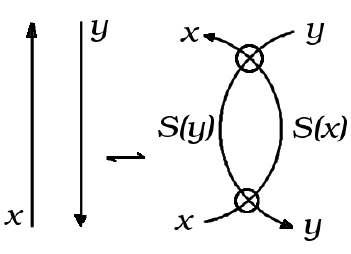}\]
Here, in the first move the contribution on the left is again 0, while on 
the right we have a contribution to the Boltzmann weight of
$v(x,y)-v(S^{-1}(y),S(x)) + v(S^{-1}(y),S(x))-v(x,y) =0$.
Similarly, in the second move
the contribution on the left is 0 while the contribution on
the right is $v(x,S(y))-v(y,S(x))+v(y,S(x))-v(x,S(y)) =0$.

\end{proof}

To deal with the virtual type III move, we now define a chain 
complex associated to the virtual biquandle operation $S:B\to B$. As in the 
classical case, we first let $C_n(B,S;\mathbb{Q})=\mathbb{Q}[B^n]$, the 
rational vector space with basis consisting of
all $n$-tuples of elements of the coloring biquandle $B$. Let
\[\partial^1_i(x_1,\dots,x_n)=(S(x_1),\dots,S(x_{i-1}),x_{i+1},\dots,x_n)\]
and
\[\partial^2_i(x_1,\dots,x_n)=(x_1,\dots,x_{i-1},S^{-1}(x_{i+1}),\dots,
S^{-1}(x_n)).\]
Let us abbreviate $\mathbf{x}=(x_1,\dots, x_n)$. Then let 
\[\partial^S_n(\mathbf{x})=\sum_{i=1}^n (-1)^i
\left(\partial^1_i(\mathbf{x})-\partial^2_i(\mathbf{x})\right)\]
and extend $\partial_n^S$ to $C_n(B,S;\mathbb{Q})$ by linearity.

\begin{proposition}\label{cc}
With the above definition, $(C_*(B,S;\mathbb{Q}), \partial^S_*)$ is a chain 
complex.
\end{proposition}

\begin{proof}
\begin{eqnarray*}
\partial^S_{n-1}(\partial^S_{n}(\mathbf{x})) & = &
\sum_{j=1}^{n-1}\sum_{i=1}^n (-1)^{i+j}
\left(\partial^1_j\partial^1_i(\mathbf{x}) -
\partial^1_j\partial^2_i(\mathbf{x}) -
\partial^2_j\partial^1_i(\mathbf{x}) +
\partial^2_j\partial^2_i(\mathbf{x}) \right) \\
& = & 
\sum_{i>j} (-1)^{i+j}
\left(\partial^1_j\partial^1_i(\mathbf{x}) -
\partial^1_j\partial^2_i(\mathbf{x}) -
\partial^2_j\partial^1_i(\mathbf{x}) +
\partial^2_j\partial^2_i(\mathbf{x}) \right) \\
& & -
\sum_{i<j} (-1)^{i+j}
\left(\partial^1_j\partial^1_i(\mathbf{x}) -
\partial^1_j\partial^2_i(\mathbf{x}) -
\partial^2_j\partial^1_i(\mathbf{x}) +
\partial^2_j\partial^2_i(\mathbf{x}) \right).
\end{eqnarray*}
Now, if $i>j$ we have
\begin{eqnarray*}
\partial^1_j\partial^1_i(\mathbf{x}) & = &
(S^2(x_1),\dots,S^2(x_{j-1}),
S(x_{j+1}),\dots,S(x_{i-1}),
x_{i+1},\dots,x_n) \\
\partial^1_j\partial^2_i(\mathbf{x}) & = &
(S(x_1),\dots,S(x_{j-1}),
x_{j+1},\dots,x_{i-1},
S^{-1}(x_{i+1}),\dots,S^{-1}(x_n)) \\
\partial^2_j\partial^1_i(\mathbf{x}) & = &
(S(x_1),\dots,S(x_{j-1}),
S^{-1}(S(x_{j+1})),\dots,S^{-1}(S(x_{i-1})),
S^{-1}(x_{i+1}),\dots,S^{-1}(x_n)) \\ 
& = &
(S(x_1),\dots,S(x_{j-1}),
x_{j+1},\dots,x_{i-1},
S^{-1}(x_{i+1}),\dots,S^{-1}(x_n)) \\
\partial^2_j\partial^2_i(\mathbf{x}) & = &
(x_1,\dots,x_{j-1},
S^{-1}(x_{j+1}),\dots,S^{-1}(x_{i-1}),
S^{-2}(x_{i+1}),\dots,S^{-2}(x_n)) \\
\end{eqnarray*}
and for $i<j$ we have
\begin{eqnarray*}
\partial^1_j\partial^1_i(\mathbf{x}) & = &
(S^2(x_1),\dots,S^2(x_{i-1}),
S(x_{i+1}),\dots,S(x_{j-1}),
x_{j+1},\dots,x_n) \\
\partial^1_j\partial^2_i(\mathbf{x})  & = &
(S(x_1),\dots,S(x_{i-1}),
S(S^{-1}(x_{i+1})),\dots,S(S^{-1}(x_{j-1})),
S^{-1}(x_{j+1}),\dots,S^{-1}(x_n)) \\
& = &
(S(x_1),\dots,S(x_{i-1}),
x_{i+1},\dots,x_{j-1},
S^{-1}(x_{j+1}),\dots,S^{-1}(x_n)) \\
\partial^2_j\partial^1_i(\mathbf{x})
& = &
(S(x_1),\dots,S(x_{i-1}),
x_{i+1},\dots,x_{j-1},
S^{-1}(x_{j+1}),\dots,S^{-1}(x_n)) \\
\partial^2_j\partial^2_i(\mathbf{x}) & = &
(x_1,\dots,x_{i-1},
S^{-1}(x_{i+1}),\dots,S^{-1}(x_{j-1}),
S^{-2}(x_{j+1}),\dots,S^{-2}(x_n)) \\
\end{eqnarray*}

Hence, the terms in the $i<j$ sum match up one-to-one with the terms
in the $i>j$ sum, but with opposite sign, and we have 
$\partial^S_{n-1}\circ\partial^S_n=0$.

\end{proof}

To find cocycles which satisfy
the condition imposed by the virtual type III move, we need to restrict our
attention to a certain subcomplex. 

\begin{definition}\textup{
Let $\rho_n:B^n\to B^n$ be given by}
\[\rho_n(x_1,\dots,x_n) = (S^{j_1}(x_1),\dots,S^{j_n}(x_n))\]
\textup{where $j_i=n-2i+1$ and let $\tau_n:B^n\to B^n$ be given by}
\[\tau_n(x_1,\dots,x_n)=(x_n,\dots,x_1).\] \textup{Let $\omega_n:B^n\to B^n$
be given by $\omega_n(\mathbf{x})=\tau_n(\rho_n(\mathbf{x}))$.}
\end{definition}

The map $\omega_n$ reverses the order of the components of its argument and
operates on each by a power of $S$ so that the sum of the power of
$S$ and the subscript of $x$ equals the position in $\omega_n$ of each 
entry. That is, the $i$th entry of $\omega_n(x_1,\dots,x_n)$ is 
$S^{2i-n-1}(x_{n-i+1})$. For example, we have
\[\omega_4(x_1,x_2,x_3,x_4) = (S^{-3}(x_4),S^{-1}(x_3),S(x_2), S^3(x_1))\]
and 
\[\omega_3(x_1,x_2,x_3)=(S^{-2}(x_3),x_2,S^2(x_1)).\]

\begin{definition}\textup{
Let $C^S_n(B,S;\mathbb{Q})$ be the subspace of $C_n(B,S;\mathbb{Q})$ generated
by elements of the form 
\[\mathbf{x}+(-1)^j\omega_n(\mathbf{x}) \quad \mathrm{where} \quad
j=\left\{\begin{array}{ll} 0 & n\equiv 0,1 \ \mathrm{mod} \ 4 \\
1 & n\equiv 2,3 \ \mathrm{mod} \ 4 \\
\end{array}\right..\] }
\end{definition}

\begin{lemma}\label{lem:sub}
For $*\in\{1,2\}$ we have $\omega_{n-1}(\partial_i^*(\mathbf{x}))
=\partial_{n-i+1}^*(\omega_n(\mathbf{x})).$
\end{lemma}

\begin{proof}
Let us denote $(y_1,\dots,y_{n-1})=\partial_i^1(x_1,\dots,x_n)$; then
we have $y_j=S(x_j)$ for $j<i$ and $y_j=x_{j+1}$ for $j>i$. 
Then
\begin{eqnarray*}
\omega_{n-1}(\partial_i^1(\mathbf{x})) & = & 
\omega_{n-1}(y_1,\dots,y_{i-1},y_i,\dots,y_{n-1}) \\
& = & \tau_{n-1}(\rho_{n-1}(y_1,\dots,y_{i-1},y_i,\dots,y_{n-1})) \\
& = & \tau_{n-1}(S^{(n-1)-1}(y_1),\dots, S^{(n-1)-2(i-1)+1}(y_{i-1}),
S^{(n-1)-2i+1}(y_i),\dots, S^{1-(n-1)}(y_{n-1})) \\
& = & \tau_{n-1}(S^{n-2}(y_1),\dots, S^{n-2i+2}(y_{i-1}),
S^{n-2i}(y_i),\dots, S^{-n}(y_{n-1})) \\
& = & \tau_{n-1}(S^{n-1}(x_1),\dots, S^{n-2i+3}(x_{i-1}),
S^{n-2i}(x_{i+1}),\dots, S^{-n}(x_n)) \\
& = & 
(S^{-n}(x_n),\dots,S^{n-2i}(x_{i+1}),S^{n-2i+3}(x_{i-1}),\dots,S^{n-1}(x_1))
\end{eqnarray*}
while
\begin{eqnarray*}
\partial_{n-i+1}^1(\omega_n(\mathbf{x})) & = &
\partial_{n-i+1}^1(S^{1-n}(x_n),\dots,S^{n-2i-1}(x_{i-1}),
S^{n-2i+1}(x_i),S^{n-2i+3}(x_{i+1}),\dots,S^{1-n}(x_1)). \\
\end{eqnarray*}
Here the entry $S^{n-2i+1}(x_i)$ is in the $n-i+1$ position, so we have
\begin{eqnarray*}
\partial_{n-i+1}^1(\omega_n(\mathbf{x})) & = &
(S^{1-n+1}(x_n),\dots,S^{n-2i-1+1}(x_{i-1}),
S^{n-2i+3}(x_{i+1}),\dots,S^{1-n}(x_1)) \\
& = & 
(S^{-n}(x_n),\dots,S^{n-2i}(x_{i+1}),S^{n-2i+3}(x_{i-1}),\dots,S^{n-1}(x_1))
\end{eqnarray*}
as required. The other case is similar.

\end{proof}

\begin{corollary}
$(C^S_n(B;\mathbb{Q}),\partial^S_n)$ is a subcomplex of 
$(C_n(B,S;\mathbb{Q}),\partial^S_n)$.
\end{corollary}

\begin{proof}
We must show that $\partial^S_n(C^S_n(B,S;\mathbb{Q}))\subset 
C^S_{n-1}(B,S;\mathbb{Q})$. For any generator
$\mathbf{z}=\mathbf{x}+(-1)^j\omega_n(\mathbf{x})\in C^S_n(B,S;\mathbb{Q}),$ 
we have
\begin{eqnarray*}
\partial^S_n(\mathbf{z}) 
& = & \partial^S_n(\mathbf{x}+(-1)^j\omega_n(\mathbf{x})) \\
& = & 
\sum_{i=1}^n (-1)^i 
\left(\partial^1_i(\mathbf{x}+(-1)^j\omega_n(\mathbf{x})) -
\partial^2_i(\mathbf{x} +(-1)^j\omega_n(\mathbf{x})) \right)\\
& = & \sum_{i=1}^n (-1)^i \left(\partial^1_i(\mathbf{x}) 
+(-1)^j\partial^1_i\omega_n(\mathbf{x}) -
\partial^2_i(\mathbf{x}) 
+(-1)^{j+1}\partial^2_i(\omega_n(\mathbf{x}))\right) \\
& = & \sum_{i=1}^n (-1)^i \left(\partial^1_i(\mathbf{x}) 
+(-1)^j\omega_{n-1}\partial^1_{n-i+1}(\mathbf{x}) -
\partial^2_i(\mathbf{x}) 
+(-1)^{j+1}\omega_{n-1}\partial^2_{n-i+1}(\mathbf{x})\right) \\
& = & \sum_{i=1}^n (-1)^i \left(\partial^1_i(\mathbf{x}) 
-\partial^2_i(\mathbf{x})\right) 
+(-1)^j\omega_{n-1}\left(\sum_{i=1}^n 
(-1)^i \left(\partial^1_{n-i+1}(\mathbf{x}) 
-\partial^2_{n-i+1}(\mathbf{x})\right)
\right) \\
& = & \sum_{i=1}^n (-1)^i \left(\partial^1_i(\mathbf{x}) 
-\partial^2_i(\mathbf{x})\right) 
+(-1)^{j+n+1}\omega_{n-1}\left(\sum_{i=1}^n 
(-1)^{n-i+1} \left(\partial^1_{n-i+1}(\mathbf{x}) 
-\partial^2_{n-i+1}(\mathbf{x})\right)
\right) \\
& = & \partial_n^S(\mathbf{x}) 
+(-1)^{j+n+1}\omega_{n-1}(\partial_n^S(\mathbf{x}))
\ \in C_{n-1}(B;\mathbb{Q}).
\end{eqnarray*}

\end{proof}

Note that 
\[\partial_n^S(\mathbf{x}\pm \omega_n \mathbf{x})
=\left\{\begin{array}{ll} 
\partial_{n-1}^S(\mathbf{x})\pm \omega_{n-1} \partial_{n-1}(\mathbf{x}) 
& n \ \mathrm{odd} \\
\partial_{n-1}^S(\mathbf{x})\mp \omega_{n-1} \partial_{n-1}(\mathbf{x}) 
& n \ \mathrm{even} \\
\end{array}\right.\]
because of the $i\leftrightarrow n-i+1$ switch, which is why we need the
$(-1)^j$ in the definition of $C_n^S(B,S;\mathbb{Q})$.

Thus, we have a chain complex $C_n^S(B,S;\mathbb{Q})$ with rational 
coefficients and  boundary maps $\partial^S_n$; passing to the dual 
gives a cochain complex $C^n_S(B,S;\mathbb{Q})$ with coboundary maps
$\delta_S^n$.

\begin{definition}
\textup{The homology and cohomology spaces}
\[H_n^S(B,S;\mathbb{Q})
=\mathrm{ker}(\partial^S_n)/\mathrm{Im}(\partial^S_{n+1}) \quad \quad
H^n_S(B,S;\mathbb{Q})=\mathrm{ker}(\delta_S^n)/\mathrm{Im}(\delta_S^{n-1})\]
\textup{of the chain complex $C_n^S(B,S;\mathbb{Q})$
are the \textit{S-homology} and \textit{S-cohomology} of the
virtual biquandle $(B,S)$.}
\end{definition}

\begin{remark}
\textup{We notice here that the $S$-chain complex is independent of the
biquandle structure on $B$; it uses $B$ only as a set and does not 
actually use the biquandle operations.}

\end{remark}

Armed with our S-homology theory, we can see what kind of function 
$v:B\times B\to \mathbb{Q}$ must be to give the same contribution to
the Boltzmann weight before and after a vIII move. 

\[\includegraphics{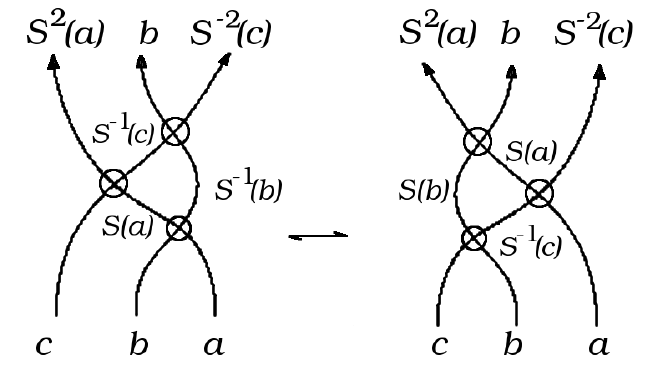}\]

Using our counting rule, the contribution to the Boltzmann weight from the 
diagram on the left is
\[c_l = v(a,b)-v(S^{-1}(b),S(a))+v(S(a),c)-v(S^{-1}(c),S^2(a)) 
+v(S^{-1}(b),S^{-1}(c))-v(S^{-2}(c),b)\]
while the contribution from the diagram on the right is
\[c_r = v(b,c)-v(S^{-1}(c),S(b))+v(a,S^{-1}(c))-v(S^{-2}(c),S(a))
+v(S(a),S(b))-v(b,S^{2}(a)).\]
Note that
\begin{eqnarray*}
\partial_3^S((a,b,c)-(S^{-2}(c),b,S^2(a))) & = &
-(b,c)+(S^{-1}(b),S^{-1}(c)) +(b,S^2(a))-(S^{-1}(b),S(a)) \\
&  & +(S(a),c)-(a,S^{-1}(c))-(S^{-1}(c),S^2(a))+(S^{-2}(c),S(a)) \\
&  & -(S(a),S(b))+(a,b) +(S^{-1}(c),S(b)) - (S^{-2}(c),b) 
\end{eqnarray*}
so that $c_l-c_r=0$ iff $v$ evaluates to $0$ on $\partial_3^S(z)$
for all $z\in C_3^S(B,S;\mathbb{Q})$. That is, the condition we need for the
contribution to the Boltzmann weight from $v$ to be invariant under $vIII$ 
moves using our counting rule is simply that $v\in C^2_S(B,S;\mathbb{Q})$.
Note that lemma \ref{lem:v} implies that the other
oriented virtual type III moves are equivalent to the one depicted here.

We are almost ready to define our virtual cocycle invariants. The final step
comes from the last oriented virtual move, the
mixed move $v$. It is here that we get interaction between the classical and 
virtual cocycles. 

Lemma \ref{lem:v} and figure 5 of \cite{CJKLS} show that in the 
presence of the oriented classical and virtual type II 
moves, a mixed move with one classical crossing is equivalent to  (1) the
same move with the classical crossing type switched and to (2) the same move 
with the orientation of the virtual strand reversed. Thus, all eight possible
oriented mixed moves are equivalent, and it is sufficient 
to consider only the one move pictured below.
 
\begin{definition}
\textup{Let $(B,S)$ be a virtual biquandle. Then a classical cocycle 
$\phi\in C^2_{YB}(B;\mathbb{Q})$ is \textit{compatible} with an $S$
cocycle $v\in C_S^2(B,S;\mathbb{Q})$ if for all $a,b,c\in B$ we have}
\begin{eqnarray*}
\phi(a,b)-\phi(S(a),S(b)) & = &
v(b,c)+v(a,S^{-1}(c))+v(S^{-1}(c),S(a^b))+v(S^{-2}(c),S(b_a)) \\
&  & -v(a^b,c)-v(b_a,S^{-1}(c)) -v(S^{-1}(c),S(b))-v(S^{-2}(c),S(a)).
\end{eqnarray*} 
\textup{If we further have 
$\phi(a,b)-\phi(S(a),S(b))=0$ for all $a,b\in B$, say that 
$\phi$ and $v$ are \textit{strongly compatible}.}
\end{definition}

Note that $v$ and $\phi$ are strongly compatible iff both are
compatible with the zero cocycle.

\begin{proposition}
If $v\in C^2_S(B,S;\mathbb{Q})$ and $\phi\in C^2_{YB}(B;\mathbb{Q})$ are
compatible, then the contributions to the Boltzmann weight of a virtual 
biquandle-colored link diagram before and after a mixed $v$ move are the same.
\end{proposition}

\begin{proof}
We simply read the condition off the picture:
\[\includegraphics{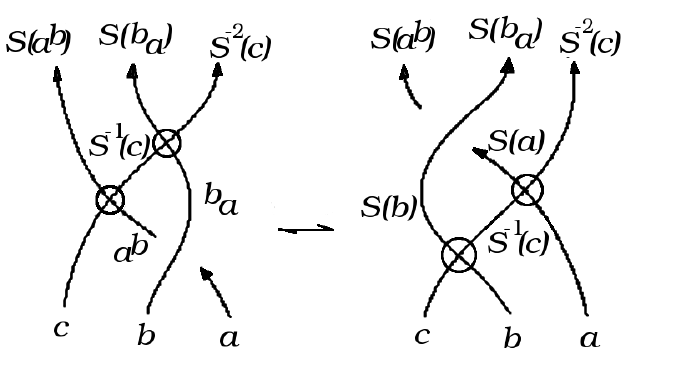}\]
\end{proof}

\begin{definition}
\textup{Let $(T,S)$ be a finite virtual biquandle, $L$ a virtual link diagram,
$\phi$ a reduced Yang-Baxter 2-cocycle of $B$ and $v$ an S 2-cocycle of
$(T,S)$ compatible with $\phi$.
For each $f\in\mathrm{Hom}(VB(L),(T,S))$ considered as a coloring
$f(L)$ of $L$, let 
\[BW_c(f)=\sum_{\mathrm{+ \ crossing}} \phi(u,o) 
-\sum_{\mathrm{- \ crossing}} \phi(u,o)
\]
and let
\[BW_v(f)=\sum_{\mathrm{virtual \ crossing}} v(r,l)-v(S^{-1}(l),S(r))
\]
where $u,o$ denote the undercrossing and overcrossing colors at classical
crossings (inbound for positive crossings, outbound for negative) and 
$r,l$ denote the right and left inbound colors at a virtual crossing.
Then the \textit{virtual Yang-Baxter 2-cocycle invariant} is the set with
multiplicities}
\[\Phi_{vYB}(L,(T,S),\phi,v)
=\{BW_c(f)+BW_v(f) \ | \ f\in\mathrm{Hom}(VB(L),(T,S))\}\]
\textup{Alternatively, we can convert the invariant into a Laurent 
polynomial (clearing denominators in $\phi$ and $v$ by scalar 
multiplication if necessary) to obtain}
\[\Phi_{vYB}=\sum_{f\in\mathrm{Hom}(VB(L),(T,S))}
t^{BW_c(f)+BW_v(f)}.\]
\textup{We will frequently omit the $L,(T,S),\phi$ and $v$ from the 
notation when these are clear from context.}
\end{definition}

If $\phi$ and $v$ are strongly compatible, the contributions from 
$v$ and $\phi$ are independent and we can sharpen the invariant to obtain a
two-variable Laurent polynomial
\[\Phi_{vYB}=\sum_{f\in\mathrm{Hom}(VB(L),(T,S)) }
t^{BW_c(f)}s^{BW_v(f)}\] or a multiset
\[\Phi_{vYB}
=\{(BW_c(f),BW_v(f)) \ | \ f\in\mathrm{Hom}(VB(L),(T,S))\}\]

By construction, $\Phi_{vYB}$ is an invariant of virtual isotopy. Curiously, 
unlike in the classical case, coboundaries can contribute nontrivially
to the Boltzmann weight for a given coloring. More specifically, as in the 
classical case, a virtual coboundary breaks down into 
weights on semiarc endpoints, i.e.
\[v(x,y)-v(S^{-1}(y),S(x))= 2f(x)+2f(y)-2f(S(x))-2f(S^{-1}(y))\]
but even if both the classical and virtual cocycles are coboundaries, 
the weights do not cancel as in the classical case unless we have
$\phi=\delta^1_{YB}(f)$ and $\frac{1}{2}v=\delta^1_{S}(f)$
for the \textit{same} $f:B\to \mathbb{Q}$. Thus, cohomologous cocycles from
either complex do not necessarily define the same virtual invariant. 

\begin{definition}
\textup{Let $C^2_D(B,S'\mathbb{Q})$ be the subspace of 
$C^2_S(B,S;\mathbb{Q})$ generated by cocycles of the form 
\[\chi_{(a,b)}+\chi_{(S^{-1}(b),S(a))}.\] Elements of $C^2_D(B,S;\mathbb{Q})$
will be called \textit{degenerate} 2-cochains.}
\end{definition}

If $v\in C^2_D(B,S;\mathbb{Q})$ then $v$ always contributes zero to the
Boltzmann weight. In particular, if two S 2-cocycles $v, v'$ differ by an 
element of $C^2_D(B,S;\mathbb{Q})\cap C^2_S(B,S;\mathbb{Q})$ then for 
every compatible classical cocycle $\phi$, the pair
$(\phi,v)$ defines the same invariant as $(\phi,v')$. 

We close
this section with a few observations. First, for every classical biquandle 
$(T,S=\mathrm{Id}_B)$, the zero virtual cocycle $v=0$ is compatible with 
every reduced Yang-Baxter cocycle $\phi$; this case is the 
classical Yang-Baxter cocycle invariant. Secondly, evaluating the polynomial 
version of the invariant at $t=1$ (or $t=s=1$ in the strong case) yields the 
virtual biquandle counting invariant $|\mathrm{Hom}(VB(L),(T,S))|$. Moreover, 
in the strong case, setting $s=t$ yields the
single-variable version of $\Phi_{vYB}$. Lastly, if $v$ and $s$ are strongly 
compatible, then the presence of any nonzero power of $s$ in the two-variable
version of the invariant indicates that the link in question is non-classical,
since a classical link will always have a total contribution of $0$ from 
virtual crossings.

\section{\large \textbf{Virtual 2-cocycle invariant examples}}\label{inv}

In this section we present some examples of virtual Yang-Baxter 
2-cocycle invariants. These were computed using \texttt{Maple} with programs
in the files \texttt{virtualbiquandles-maple.txt}, 
\texttt{biquandles-maple.txt}, and \texttt{yangbaxtercohomology.txt},
available at \texttt{www.esotericka.org/quandles}.

Our \texttt{Maple} notation uses Gaussian integers for crossing labels
in virtual signed Gauss codes with imaginary part 0 for positive crossings, 
imaginary part $\pm i$ for negative crossings and imaginary part $\pm 2i$ for 
virtual crossings. The absolute value of real part of the label gives the
crossing label; we use negative entries for undercrossings and right-hand 
virtual labels, and positive entries for overcrossings and left-hand virtual 
labels. A zero entry separates components. For example, the
virtual signed Gauss code $U1^-R2\ O3^+O1^-L2\ U3^+$ in our \texttt{Maple} 
notation becomes 
\begin{center}\texttt{[-1-I,-2-2*I,3,1+I,2+2*I,-3,0].}\end{center}

The program \texttt{vblist} takes a biquandle matrix (actually, a 4-tuple of
operation matrices) and computes the automorphism group, then outputs one 
representative of each conjugacy class, with a permutation $\sigma\in S_n$
represented by the vector $[\sigma(1),\dots, \sigma(n)]$.

The program \texttt{vbhomlist} takes a virtual Gauss code, a biquandle matrix
and a permutation \texttt{S} and computes the set of all colorings of the 
virtual link represented by a virtual Gauss code by the virtual biquandle. 
\texttt{vcocycles} finds a basis for the space of $S$ 2-cocycles
given a permutation $S$. These  2-cocycles are represented by 
$n^2$-component vectors where the $n(i-1)+j$ entry represents the
coefficient of $\chi_{(ij)}$, where $n$ is the cardinality of $B$.  
\texttt{ybcompatfind} takes a biquandle matrix \texttt{A} and a permutation
\texttt{S} and finds a list of pairs of compatible classical and 
non-degenerate $S$-cocycles. Finally, \texttt{vbinv} takes a virtual signed 
Gauss code \texttt{g}, a biquandle matrix \texttt{A}, a permutation 
\texttt{S}, a classical cocycle vector $\texttt{phi}$ and an $S$-cocycle 
vector \texttt{v} and computes the polynomial virtual Yang-Baxter invariant.
The program checks for strong compatibility and automatically computes
the two-variable invariant (with variables \texttt{s,t}) if the cocycles are 
strongly compatible, or the one-variable invariant (with variable \texttt{T})
if they are only weakly compatible.

\begin{example}
\textup{The trivial biquandle of cardinality $n$, $T_n=\{1,2,\dots,n\}$
with $a^b=a^{\bar{b}}=a_b=a_{\bar{b}}=a$ has reduced Yang-Baxter
cocycles spanned by $\chi_{(ij)}$ where $i\ne j$. If $K$ is a 
single-component virtual link (i.e., a virtual knot), then the only classical 
biquandle colorings by $T_n$ are
the constant colorings and the classical 2-cocycle invariant is equal
to the counting invariant. However, if we endow
$T_n$ with a nontrivial automorphism $S$ (and any permutation will do since 
$\mathrm{Aut}(T_n)=S_n$), then a compatible classical cocycle
can contribute non-trivially to the virtual Yang-Baxter cocycle invariant.
Here the virtual trefoil has $\Phi_{vYB}(K)=3t$ with respect to the
Yang-Baxter cocycle $\phi=\chi_{(13)}+\chi_{(21)}+\chi_{(23)}$ and the
zero $S$-cocycle $v=0$ of $(T_3,(123))$.}
\[\includegraphics{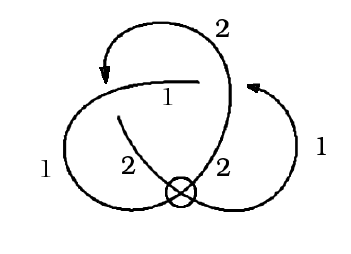}\quad
\includegraphics{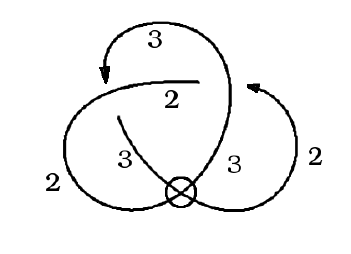}\quad
\includegraphics{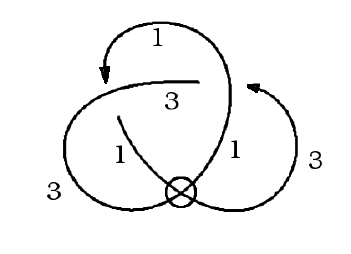}\quad
 \] 
\end{example}

\begin{example}
\textup{The virtual biquandle $T=\{1,2,3\}$, $S=(23)$ with biquandle 
operation matrix below has S-cocycle $v=\chi_{(12)}+\chi_{(13)}$ strongly 
compatible with the zero classical cocycle $\phi=0$.
The virtual link $L$ below has virtual Yang-Baxter cocycle invariant 
$\Phi_{vYB}=2s^{-2}+5+2s^2$ with respect to this $\phi$ and $v$. Thus, 
this invariant detects the non-classicality of $L.$ We list the colorings 
in the table.}
\[M_T=\left[\begin{array}{ccc|ccc}
1 & 1 & 1 & 1 & 1 & 1 \\
2 & 3 & 3 & 2 & 3 & 3 \\
3 & 2 & 2 & 3 & 2 & 2 \\ \hline
1 & 1 & 1 & 1 & 1 & 1 \\
3 & 3 & 3 & 3 & 3 & 3 \\
2 & 2 & 2 & 2 & 2 & 2 
\end{array}\right] \quad
\raisebox{-0.5in}{
\includegraphics{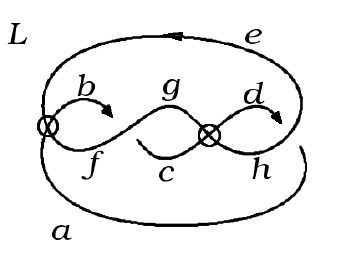}} \quad
\begin{array}{|cccccccc|}
\hline
a & b & c & d & e & f & g & h \\ \hline
1 & 1 & 1 & 1 & 1 & 1 & 1 & 1 \\
1 & 1 & 1 & 1 & 2 & 3 & 2 & 3 \\
1 & 1 & 1 & 1 & 3 & 2 & 3 & 2 \\
2 & 3 & 3 & 2 & 1 & 1 & 1 & 1 \\
2 & 3 & 2 & 3 & 2 & 3 & 2 & 3 \\
2 & 3 & 2 & 3 & 3 & 2 & 3 & 2 \\
3 & 2 & 2 & 3 & 1 & 1 & 1 & 1 \\
3 & 2 & 3 & 2 & 2 & 3 & 2 & 3 \\
3 & 2 & 3 & 2 & 3 & 2 & 3 & 2 \\ \hline
\end{array}\]
\end{example}

\begin{example}
\textup{Let $T=\mathbb{Z}_3$ with $a^b=a^{\bar{b}}=a_b=a_{\bar{b}}=2a$ mod 3.
Then $T$ is an Alexander biquandle with $s=t=2$, and 
$\mathrm{Aut}(T)=\{(1),(12)\}$. $T$ has Yang-Baxter cocycle 
$\phi=\chi_{(13)}+\chi_{(23)}$ which is also an $S$-cocycle for the 
automorphism $S=(12)\in \mathrm{Aut}(T)$. Moreover, $\phi$ and $v=\phi$ are 
strongly 
compatible. The virtual link below has virtual Yang-Baxter invariant 
$\Phi_{vYB}=5+2s^{-2}+ 2s^2t^{-2}$ with respect to $\phi=v$.}
\[\includegraphics{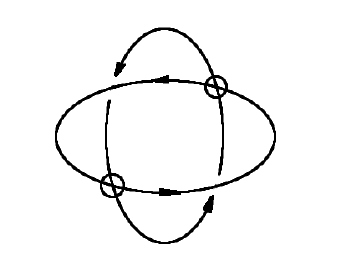}\]
\end{example}

\begin{example}
\textup{The Alexander quandle $T=\mathbb{Z}_6=\{1,2,3,4,5,6\}$ with}
\[a^b=a^{\bar{b}}=5a+2b \ \mathrm{mod} \ 6,\quad 
a_b=a_{\bar{b}}=a \ \mathrm{mod} \ 6\]
\textup{with $S=(26)(35)$ has an $S$-cocycle} 
\[v=-\chi_{(14)}-2\chi_{(16)}+2\chi_{(22)}+2\chi_{(23)}+\chi_{(24)}
+2\chi_{(25)}-\chi_{(34)}+\chi_{(42)}+\chi_{(43)}
\] 
\textup{weakly compatible with the Yang-Baxter cocycle}
\[\phi=\chi_{(21)}+\chi_{(23)}+\chi_{(24)}+\chi_{(26)}
-\chi_{(43)}-\chi_{(46)}-\chi_{(61)}-\chi_{(64)}\]
\textup{In fact, $\phi$ is a coboundary, but the virtual cocycle invariant 
defined by $\phi$ and $v$ includes nonzero contributions from $\phi$. For
example, the virtual Hopf link below has invariant value 
$\Phi_{vYB}=3T^{-1}+6+3T$;
four of the twelve colorings have nonzero Boltzmann weight contributions
from $\phi$.}
\[\raisebox{-0.5in}{\includegraphics{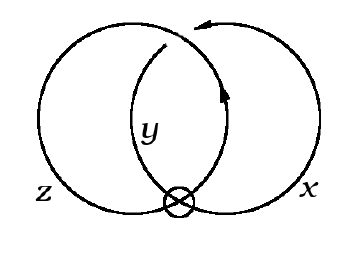}} \quad
\begin{array}{|ccc|} \hline 
x & y & z \\ \hline
1 & 1 & 1 \\
2 & 6 & 1 \\
3 & 5 & 1 \\ 
4 & 4 & 1 \\
5 & 3 & 1 \\
6 & 2 & 1 \\\hline
\end{array}\ 
\begin{array}{|ccc|} 
\hline x & y & z \\ \hline
1 & 1 & 4 \\
2 & 6 & 4 \\
3 & 5 & 4 \\
4 & 4 & 4 \\
5 & 3 & 4 \\
6 & 2 & 4 \\ \hline
\end{array}
\]
\end{example}

\section{\large \textbf{Questions for future research}}

We have only scratched the surface of the subject of virtual Yang-Baxter 
cocycle invariants. In this section we collect a few questions for future
research.

As is well known, quandle cocycle invariants can be defined not only for
knotted circles in $S^3$ and other 3-manifolds, but also for knotted 
surfaces in $\mathbb{R}^4$
or other 4-manifolds. The surface knot equivalents of virtual knots are
known as \textit{abstract surface knots}, equivalence classes of diagrams
obtained by gluing together boxes containing triple points, pairs of crossed 
sheets, etc., under the equivalence relation defined by the Roseman moves. 
See \cite{KK}.
What modifications to our scheme are necessary to define virtual cocycle
invariants for abstract surface knots?

The fact that we have two separate homology theories contributing to
our invariants seems unsatisfying, especially as the compatibility condition 
is not expressed in terms of either coboundary map. Moreover, the fact that
cohomologous cocycles do not define the same invariant seems troubling, as 
does the fact that the degeneracy condition for $S$ 2-cocycles is not just
the coboundary condition.
Perhaps there is a deeper cohomology theory which unifies both the classical
Yang-Baxter and $S$-homology theories? It is because we suspect that such a
theory awaits uncovering that we chose the term ``$S$-homology'', preferring 
to reserve the term ``virtual biquandle homology'' for this deeper theory.

Other counting rules for contributions to the Boltzmann weight from both
classical and virtual crossings are possible, of course, and choosing 
different rules may result in different theories, defining potentially
different invariants. We chose our rule for counting contributions from 
virtual crossings because it handles the $vI$ and $vII$ moves easily; 
different rules might have other advantages and disadvantages. It is 
worth noting that the obvious alternative counting rules appear to suffer 
from the same deficiencies as the rule we ultimately chose, and most have the 
added disadvantage of requiring quotients to handle $vI$ and $vII$ moves.
For example, if we count just $v(r,l)$ instead of $v(r,l)-v(S^{-1}(l),S(r))$
we can use the same boundary map $\partial_n^S$, but two levels of quotients 
are required to handle the $vI$ and $vII$ moves.

The virtual biquandle definition we used is the simpler of the two given in
\cite{KM}; instead of using the unary operations $S$ and $S^{-1}$ at virtual 
crossings, one can define a pair of binary operations at a virtual crossing, 
so that the virtual biquandle structure consists of six binary operations. 
A theory of cocycle invariants associated to this presumably richer virtual 
biquandle structure should prove interesting as well.

\end{document}